\documentclass[12pt]{article}
\usepackage{amsfonts,amssymb,amsmath,amscd}

\def\be{\begin{equation}}
\def\ee{\end{equation}}
\def\bea{\begin{eqnarray}}
\def\eea{\end{eqnarray}}

\def\C{{\mathbb C}}

\def\hbar{{\xi}}

\begin{document}

\thispagestyle{empty}

\hfill \today 

\vspace{2.5cm}

\begin{center}
\sc{\LARGE    QUANTUM BASES IN $U_q(g)$ }
\end{center}

\bigskip\bigskip

\begin{center}
Enrico Celeghini
\end{center}

\begin{center}
{\sl Departimento di Fisica, Universit\`a  di Firenze and
INFN--Sezione di
Firenze \\
I50019 Sesto Fiorentino,  Firenze, Italy}\\
\medskip

{e-mail: celeghini@fi.infn.it}

\end{center}

\bigskip

\bigskip

\begin{abstract}

This paper is devoted to analyze, inside the $\infty$-many possible
bases of a Quantum Universal Enveloping Algebra $U_q(g)$, those
that can be considered as ``more equal then others'', like orthonormal 
bases in the Euclidean spaces.

The only possible element of selection has been found to be a privileged
connection with the corresponding bialgebra.
A new parameter $z' \in \C$ --independent from the 
$z:= {\rm log}(q) \in \C$ that
defines $U_q(g)$-- has thus been introduced. 
Each value of such new parameter $z'$ defines one of these bases (we call 
quantum bases) and determines, independently from $z$, its 
commutation relations.
Both $z$ and $z'$ are, on the contrary, necessary to fix 
the relations between the basic set and its co-products. We have thus --for 
each pair $\{z,z'\}$-- one $n$-dimensional quantum basis 
$(g_{z'},\Delta_{zz'})$, that describes $U_q(g)$.

Three cases are 
particularly relevant: the analytical basis $g_{z}$, where $z'=z$,
the Lie basis $g_{0}$ obtained for $z'=0$ (where both the basic set and 
its co-product close Lie-like commutation relations and the non primitive
co-product describe an interaction) and the canonical/crystal basis 
$g_{\infty}$, limit for 
$z' \to \infty$ in the Riemann plane of the generic quantum basis.

To simplify the exposition, we discuss in details the easily generalizable
case of $U_q(su(2))$.

\end{abstract}

\vskip 1cm


\vskip 0.4cm


\vfill\eject

\section{Introduction}

Quantum Universal Enveloping Algebras $U_q(g)$ are $\infty$-dimensional
objects, not easy to manage. More, there is a large freedom in their
description that can be realized by means of
$\infty$-many equally good basic sets of $n$ elements.
It is thus difficult also to determine if two $U_q(g)$ are 
equivalent or not. To simplify the task we need to individuate,
inside the possible bases --i.e. inside the $n$-dimensional sets
that allow to build as a set of monomials of ordered powers
a Poincar\'e-Birkhoff-Witt (PBW) basis\cite{CP}--  same privileged ones 
with peculiar properties. 
It is a little like as to look for orthonormal bases in the Euclidean spaces .

For the Lie case the situation is easy:
almost everything can be done in applications without
taking into account the Universal Enveloping Algebras $U(g)$, 
because the space of generators is clearly defined as the unique
basic set that is
closed under commutation relations and
additive in the direct product space.
$U(g)$ can thus be disregarded for practical 
applications, but the situation is quite different in quantum algebras.
The set of possible equivalent bases is indeed quite
larger and no bases exists that are primitive nor
closed under commutation relations.
Linear transformations loose thus their privileges and 
the requirement that the relations between two bases are
invertible is the unique that survives. 

While in mathematics the interest is focused on $U_q(g)$ and not on its
bases, in applications it is essential to 
individuate, inside $U_q(g)$, the well precise elements to be connected
to the physical observables as --for instance-- in $U(su(2))\,$, where the Lie
algebra generators correspond to the angular momentum operators.

Our starting point is the one-to-one connection between $U_q(g)$ 
and the related bialgebra \cite{CBO}.
The $z$ parameter of the bialgebra determines the whole $U_q(g)$, 
but the situation is different when we are interested not only in the
$\infty$-dimensional $U_q(g)$ but also in a 
$n$-dimensional basic set inside it.
  
The problem is well known and a basis exists in quantum algebras that is
called canonical as it is suggested
\cite{Lu,Ka} as a possible
basis that emerges among all the others. However this basis
(called also crystal basis) is regular inside the representations, 
but it is highly discontinuous near the highest weights.

Note also that a choice among the bases has been informally done by the
scientific community: without any official justification almost always
$U_q(su(2))$ is written in a basis similar to that described in 
eqs.(\ref{coprodJ}, \ref{commJ}), related to the analyticity of 
the whole scheme.
Analyticity indeed allows to select this basis, we call 
analytical basis. It is built by means of a 
perturbative approach 
introducing order by order
the modifications required by the consistency of the bialgebra with
the Hopf algebras postulates \cite{CBO}.

However analytical basis does not look the only possible choice
as canonical/crystal basis \cite{Kang, MS}
and, in particular, Lie basis
have properties that look interesting for applications.
What is nevertheless clear is that, among the plethora of invertible 
functions of the Lie generators a subset clearly appears: 
the basic sets related to the bialgebra. We call these sets quantum bases
and we suggest that --as
they are derived analytically from the co-commutator $\delta$-- 
they are ``more equal then others''.

Each of these quantum bases is characterized by a value of a new parameter 
$z' \in \C$ totally independent from the $z:= {\rm log}(q) \in \C$ 
that defines the Hopf algebra $U_q(g)$.
Lie, canonical/crystal and analytical bases correspond to
$z'=0$, $z'=\infty$ and $z'=z$, respectively.

Analyticity is assumed in both these parameters $z$ and $z'$.
It is indeed obvious that the parameter $q=e^z$, normally used in quantum
groups,
is not the best analytical candidate:
there is in quantum algebras a symmetry under the interchange 
$q \leftrightarrow q^{-1}$
that is not well described in this parametrization and the
Lie-Hopf algebra is obtained for $Im(q)= 0$, $Re(q)=1$ in the complex 
plane $q$ that doesn't look a sufficiently symmetric 
point for a particular so relevant case. 
The situation is quite better with $z$ and $z'$.
Indeed everything in quantum algebras is written in terms of 
polynomials and exponentials in $z$ and $z'$, all
entire functions: the
structure is thus always perfectly defined for $z$ and $z'$ finite while all  
singularities are at the infinity in the Riemann sphere.
Everything is invariant under $z' \to -z'$ and
a generalized co-commutativity (i.e. an 
invariance under both $z \leftrightarrow -z$ and interchange in 
the order of the spaces in the co-products) can be introduced \cite{BCO}.
Finally Lie basis is written for $z'=0$ and the canonical/crystal
basis is defined in the point conjugate 
~i.e. the point at infinity. 
A sort of conjugation is in this way prefigured between them.

Inside $U_q(g)$ --i.e. for a fixed value of $z$--
four cases are considered: 

1) a basis $g_z$ where commutation relations are written in 
terms of the same parameter $z$ that defines $U_q(g)$, i.e. such that
the parameter 
$z'$ that define the commutation relations coincides with $z$. 
It is the basis we called analytical basis and also
the ``traditional'' basis of $U_q(g)$. 

2) a set of quantum bases $g_{z'}$ with $z' \neq z$ finite,
where commutation relations are written in terms of this new parameter
while co-products depend from both $z$ and $z'$.

3) one basis $g_{0}$ where the basic set is closed under 
commutation relations, i.e. it exhibits a Lie-like behavior.  The parameter 
$z'={\rm log}(q')$ introduced in the commutation relations is thus $z'=0$.
The co-product, on the contrary, is determined
by the coalgebra to be not primitive.
This basis could describe a conserved Lie symmetry where the presence of 
$z\neq 0$ is a signal of one interaction in the
composed system \cite{CO}.

4) a basis $g_{\infty}$ obtained in the limit $z'\to\infty$, i.e.
the canonical/crystal basis.

All these are quantum bases of the same Hopf algebra $U_q(g)$ and, one of 
them known, it is a simple
problem of change of basis --obtained in the following from 
representations theory-- to have a complete description of $U_q(g)$ 
in each of them.

As the essential structures of a quantum algebra are the commutation relations
and the co-products, we do not discuss here the other properties of 
Hopf algebras, co-unity and antipode, because the procedure 
is similar and, up to now, these objects have not played a role in 
applications \cite{CP}.  

We restrict ourselves in this paper to ``the basic 
example'',
standard deformation of $U(su(2))$ for two reasons.
First, because more then half of the papers in the applications of the 
Lie algebras
are related to $su(2)$ --both as $su(2)$ is elementary but
not so much to be trivial and because the invariance under 3-dimensional
rotations is the most fundamental one in physics-- and 
in quantum algebras the situation is not different as
the fundamental representation of $SU_q(2)$, introduced in St. Petersburg 
to describe the inverse scattering \cite{FRT},
could be considered the starting
point of all this field of research.
Second, a general discussion --with a set of theorems
for any $U_q(g)$ and any (multidimensional) bialgebra-- should be
quite heavy without adding nothing to the understanding of the problem and of
the proposed solution.

\section{Analytical basis $su_{z}(2)$ of $U_q(su(2))$}

Let us begin fixing the notations for the Lie algebra $su(2)$ (that,
in our notation, would be written $su_{0}(2)$): 

\be\label{commL}
[ L_3, L_\pm] = \pm L_\pm \qquad
[L_+,L_-] ~= 2 L_3 \,,
\ee

\begin{equation}\label{reprL}
\begin{array}{l}
L_3 |j,m\rangle = |j,m\rangle\, m\\[0.3cm]
L_{\pm}|j,m\rangle = |j,m \pm 1\rangle \;\sqrt{(j+1/2)^{2} - 
(m\pm 1/2)^{2}} 
\end{array}\ee
where the parameter $j$, related to the Casimir invariant $L^2$,
belongs to $U(su(2))$ as it
can be written as a series in its PBW basis:
\begin{equation}\label{CasimirL}
j = \frac{1}{2} \sqrt{1+\frac{1}{2} \{L_+,L_-\}+ L_3^2} - \frac{1}{2}\;.
\end{equation}

\noindent The algebra $su(2)$ is usually implemented to be the 
Lie-Hopf algebra $(su_0(2), \Delta_{(0)})$ 
introducing the primitive co-product $\Delta_{(0)}$
that describes the additivite composed systems:
\be\label{DL}
\Delta_{(0)}(L_3):= \,L_3 \otimes 1 + 1 \otimes L_3\qquad
\Delta_{(0)}(L_\pm):= \,L_\pm \otimes 1 + 1 \otimes L_\pm \;
\ee 
\noindent and is homomorphic to the algebra as
\[
[ \Delta_{(0)}(L_3), \Delta_{(0)}(L_\pm)] = \,\pm \;\Delta_{(0)}(L_\pm) \qquad
[\Delta_{(0)}(L_+),\Delta_{(0)}(L_-)] = \,2 \;\Delta_{(0)}(L_3) \,.
\]

Now, starting from the Lie algebra $su(2)$, the standard bialgebra 
$(su(2),\delta)$ can 
be constructed (see e.g. \cite{CP}) including in the game
the standard co-commutators $\delta$:
\be\label{dL}
\delta(L_3) = 0,\qquad \delta(L_\pm) = z \,(L_3 \otimes L_\pm -
L_\pm \otimes L_3) \;.
\ee
\noindent These co-commutators are defined in such a way that 
$\{\Delta_{(0)}(J_m) + 
\,\delta(J_m)\}$ for 
$(m=3,+,-)$ are homomorphic to (\ref{commL}) up to 
third order in $\{J_m \otimes 1\}$ and $\{1 \otimes J_m\}$:
\begin{equation*}
\begin{split}
[& \Delta_{(0)} (L_3)+\,\delta(L_3)\,,\; \Delta_{(0)}(L_\pm) 
+ \,\delta (L_\pm)] 
= \,\pm \,\left(\Delta_{(0)}(L_\pm) + \,\delta(L_\pm)\right) + o(3  ) \\[0.3cm]
[&\Delta_{(0)}(L_+) + \,\delta(L_+) \,,\;\Delta_{(0)}(L_-) +\,\delta(L_-)] 
\,= \,2 \,\left(\Delta_{(0)}(L_3) + \,\delta(L_3)\right)+o(3  )\,.
\end{split}
\end{equation*}
In other words the commutation relations of the composed system
are invariant up to the second order in the generators under
the addition of $\delta$ to $\Delta_{(0)}$.
The (\ref{dL}) describe the behavior of the standard Hopf algebra 
$U_q(su(2))$ near the Lie limit and, because of the Hopf algebras postulates, 
are sufficient to determine entirely it.
However, while the bialgebra is a 3-dimensional object,
the Hopf algebra is an $\infty$-dimensional one and it can be described 
by means of $\infty$-many arbitrary basic sets.
The problem in applications is to individuate ``the good one''
i.e. the basic set to relate to the physical observables as it happened
in $U_q(su(2))$ where the inverse scattering
requires a well precise basis (that 
is just the basis obtained from the bialgebra by analytical 
continuation).
  
Starting from $(g, \delta)$ --as shown in \cite{CBO}-- this analytical basis
has been obtained introducing order by order only 
the changes 
required by the Hopf algebra postulates without free parameters whatever. 

We assume that its commutators can be written 
as a formal series in $\{J_m\}$ and its
co-products as formal series in $\{1 \otimes J_m\}$ and $\{J_m \otimes 1\}$
and from the co-commutators of the Lie bialgebra the analytical
co-products are
reconstructed order by order introducing only the contributions imposed by
the co-associativity
\be\label{perturbcoass}
\sum_{j=0}^k 
(\Delta_{(j)} \otimes 1 - 1 \otimes \Delta_{(j)}) \circ
\Delta_{(k-j)} = 0 \qquad \forall k \;\,,
\ee
where the co-products are
\[
\Delta = \sum_{k=0}^\infty \Delta_{(k)} = \Delta_{(0)} + \Delta_{(1)} + 
\Delta_{(2)}\, \dots
\]
with $\Delta_{(k)}$ polynomial of order $k+1$ in $\{J_m \otimes 1\}$ 
and $\{1 \otimes J_m\}$.

The $\Delta_{(k)}$ thus obtained are then summed to:
\be\label{coprodJ}
\begin{array}{l}
\Delta(J_3) = 1 \otimes J_3 + J_3 \otimes 1 \\[.3cm]
\Delta(J_\pm) = q^{J_3} \otimes J_\pm + J_\pm \otimes q^{-J_3} \,.
\end{array}
\ee

Finally the homomorphism of the co-product, always
perturbatively written as
\[
\Delta_{(k)} ([J_m, J_n]) = \sum_{j=0}^k [ \Delta_{(j)}(J_m), 
\Delta_{(k-j)}(J_n)]
\qquad \forall k \quad \,, 
\]
has been used to derive order by order the analytical commutation relations 
always disregarding
at each order all possible modifications independent from $\delta$ 
as related to a non linear basis inside the Lie limit $U(su(2))$.
Summing the series we have:
\be\label{commJ}
[J_3,J_{\pm}] = \pm J_{\pm} \qquad\quad
[J_+,J_-] =  ~[2J_3]_q
\ee
where $[n]_q$ is related to the usual $q$-number \cite{GR} but does not
coincide with it:
\be\label{qnumber}
[n]_q := \frac{{\rm sinh}(z n) }{z}\; .
\ee
The analytical approach indeed implies the absence in the denominator of 
$[n]_q$ of the hyperbolic sine that would disconnect in eqs.(\ref{commJ})
the powers of $z$ from those of the $J_3$~.

Note that in the Lie case the structure constants are in mathematics assumed
adimensional together with the generators but this is 
not what we find from physics.
Physics, indeed, prescribes for simple algebras that both the structure 
constants and the generators are dimensional and have the same dimension. 
Only with an arbitrary change of normalization all objects can be reduced
to be adimensional. 
As we require that the physical observables can be interpreted as
generators also after the
quantization, eq.(\ref{qnumber}) is consistent and requires 
${\rm dim}[z] ={\rm dim}[J_3]^{-1}$, while the usual $q$-number cannot be
accepted. 

As shown by (\ref{commJ}), in the analytical basis the $z$ introduced 
by the bialgebra is the only one parameter that appears. It describes indeed
not only $U_q(su(2))$ but also the commutation relations of the basis.

To look for other bases we now consider the representations,
as the space of the representations does not change in the quantization.
The same analytical/perturbative procedure can be applied to
eqs.(\ref{reprL}) to obtain
\begin{equation}\label{reprJ}
\begin{array}{l}
J_3 |j,m\rangle = |j,m\rangle \,m\\[0.3cm]
J_{\pm}|j,m\rangle = |j,m \pm 1\rangle 
\;\sqrt{ [j+1/2]_q^{\,\;2} ~-~ 
[m\pm 1/2]_q^{\,\;2}}
\;/\;{\sqrt{[1]_q}} 
\end{array}\ee
where $j$ can be alternatively written in the PBW basis of $su_{z}(2)$.

\section{Quantum Bases $su_{z'}(2)$}

Let us now consider a different Quantum Universal Enveloping Algebra
$U_{q'}(su(2))$,
defined in terms of $\{K_3,K_\pm\}$, 
related to the same bialgebra (\ref{commL},\ref{dL}) but corresponding to 
a different 
quantization parameter $z'$ (or, equivalently, $q':= e^{z'}$).
Obviously equations (\ref{commJ}) and (\ref{reprJ})
remain valid with the change $q \to q'$:

\begin{equation}\label{commK}
\begin{array}{l}
[K_3,K_{\pm}] = \pm K_{\pm} \\[0.3cm]
[K_+,K_-] =  ~[2K_3]_{q'}
\end{array}\end{equation}

\begin{equation}\label{reprK}
\begin{array}{l}
K_3 |j,m\rangle = |j,m\rangle \,m\\[0.3cm]

K_{\pm}|j,m\rangle = |j,m \pm 1\rangle ~~\sqrt{[j+1/2]_{q'}^{\,2} - 
[m\pm 1/2]_{q'}^{\,2}}
\;/\;{\sqrt{[1]_{q'}}} \;.
\end{array}\end{equation}

Also the co-products of $\{K_m\}$ look like the ones in eq. 
(\ref{coprodJ}) with $q \to q'$ when the Hopf algebra $U_{q'}(su(2))$
is considered. However, while eqs.(\ref{commK},\ref{reprK}) remain valid,
$\{\Delta(K_m)\}$ appear completely different if we adopt the different point 
of view and we describe by means of the $\{K_m\}$ our original
Hopf algebra $U_q(su(2))$.
Indeed --as they are defined on the same representations space--
the $\{K_i\}$ define a PBW basis for $U_q(su(2))$ also.
Comparing eqs. (\ref{reprJ}) and (\ref{reprK}) we can rewrite 
$\{K_m\}$ in terms of $\{J_m\}$ as 
\begin{equation}\label{mapKtoJ}
K_3 = J_3 \;\qquad\quad
K_{\pm} =\sqrt{\frac{[1]_{q}}{[1]_{q'}}} \;\;
\sqrt{\frac{
[j+1/2]_{q'}^{\,2} - [J_3\mp 1/2]_{q'}^{\,2}}
{[j+1/2]_q^{\,2} - [J_3\mp 1/2]_q^{\,2}}}\;\;
J_\pm \;\;.
\end{equation}

Because we have assumed to be inside $U_q(su(2))$, we have 
from (\ref{mapKtoJ})
\be
\Delta(K_3) = \Delta(J_3)
\ee
\begin{equation}\label{mapDKtoDJ}
\Delta(K_{\pm}) = \sqrt{\frac{[1]_{q}}{[1]_{q'}}} \;\;
\sqrt{\frac{
\Delta([j+1/2]_{q'}^{\,2}) - \Delta([J_3\mp 1/2]_{q'}^{\,2})}
{\Delta([j+1/2]_q^{\,2}) - \Delta([J_3\mp 1/2]_q^{\,2})}}\;
\Delta(J_{\pm})
\end{equation}
where all the {\it rhs} is defined in terms of the $q$-depending 
$\{\Delta(J_m)\}$ of eqs. (\ref{coprodJ}) and the $\{\Delta(K_m)\}$ are
written in function of the $\{\Delta(J_m)\}$. 

But relations (\ref{mapKtoJ}) are easily invertible
\begin{equation}\label{mapJtoK}
J_3 = K_3 \qquad
J_{\pm} = \sqrt{\frac{[1]_{q'}}{[1]_{q}}}\;\;
\sqrt{\frac{
[j+1/2]_{q}^{\,2} - [K_3\mp 1/2]_{q}^{\,2}}
{[j+1/2]_{q'}^{\,2} - [K_3\mp 1/2]_{q'}^{\,2}}}\;\;
K_{\pm}
\end{equation}
and -- as the $\Delta(J_m)$ are written in terms of $\{J_m\otimes 1\}$ and
$\{1 \otimes J_m\}$--
it is sufficient to substitute in the {\it rhs} of 
eq (\ref{mapDKtoDJ}), $\{J_m\otimes 1\}$ and $\{1\otimes J_m\}$ with 
their expressions, as given by eq.(\ref{mapJtoK}), to obtain 
$\{\Delta(K_m)\}$ in terms of $\{K_m \otimes 1\}$ and $\{1\otimes K_m\}$.

To be more clear, let us consider the first elements of the series 
in (\ref{mapKtoJ}):

\[
K_\pm = J_\pm + 
\frac{{z'}^2-z^2}{3} {\cal S}(J_3^2 J_\pm) + 
\frac{{z'}^2-z^2}{6} {\cal S}(J_\pm^2 J_\mp) + o(4) . 
\]
where $o(4)$ means contributions to the fourth or higher order in $\{J_m\}$ and
${\cal S}$ is the symmetrizer i.e. a linear operator 
that, for all $n$, operates as
\[
{\cal S}(O_1 O_2 \dots O_n):=  \frac{1}{n!} \sum_{\sigma_p \in S_n} 
\sigma_p \,(O_1 O_2 \dots O_n)
\]
with $S_n$ permutation group on $n$ elements.

The co-products of $\{K_m\}$ are then rewritten in terms of 
the $\{\Delta(J_m)\}$ of (\ref{coprodJ}) as
\[
\Delta(K_\pm) = \Delta(J_\pm)  + 
\frac{{z'}^2-z^2}{3} \Delta_{(0)}\bigl({\cal S}(J_3^2 J_\pm)\bigr)+ 
\frac{{z'}^2-z^2}{6} \Delta_{(0)}\bigl({\cal S}(J_\pm^2 J_\mp)   
\bigr) + o(4) \;.
\]
Now we can substitute in the {\it rhs} 
to the $\{J_m\otimes 1\}$ and $\{1\otimes J_m\}$ their expressions in 
function of $\{K_m\otimes 1\}$ and $\{1\otimes K_m\}$ respectively, as
given by eq. (\ref{mapJtoK}), and  we obtain

\begin{equation}\label{perturbDK2}
\begin{array}{l}
\Delta(K_3) = 1 \otimes K_3 + K_3 \otimes 1 \\[0.3cm]
\Delta(K_\pm) = q^{K_3} \otimes K_\pm + 
K_\pm \otimes q^{-K_3}  \;+  \\[0.3cm]

\frac{{z'}^2-z^2}{3} \Bigl[\Delta_{(0)}({\cal S}(K_3^2 K_\pm))- 
1 \otimes {\cal S}(K_3^2 K_\pm)- {\cal S}(K_3^2 K_\pm)\otimes 1 
\Bigr] \;+ \\[0.3cm]

\frac{{z'}^2-z^2}{6} \Bigl[\Delta_{(0)}({\cal S}(K_\pm^2 K_\mp))- 
1 \otimes {\cal S}(K_\pm^2 L_\mp)- {\cal S}(K_\pm^2 K_\mp)\otimes 1 
\Bigr]
 + o(4) \,. 
\end{array}\end{equation}

Note that from (\ref{perturbDK2})
\begin{equation}\label{perturbDK3}
\begin{array}{l}
\Delta(K_3) = 1 \otimes K_3 + K_3 \otimes 1 + o(3)\\[0.3cm]
\Delta(K_\pm) = (1 \otimes K_\pm + K_\pm \otimes 1) +\, 
z\, (K_3 \otimes K_\pm - K_\pm \otimes K_3) + o(3)
\end{array}\end{equation}
\noindent that show that $\{K_m\}$, as required, define a basic set 
of $U_q(su(2))$ and are not related to $U_{q'}(su(2))$.

\section{Lie basis $su_{0}(2)$}

We call Lie basis $su_0(2)$ of $U_q(su(2))$ the basis where the 
commutation relations are Lie-like
both for the basic elements and for their co-products: all the problems of
consistency with the Hopf algebra has been moved to the form of the co-product.
It is the special case of quantum bases for $z'=0$.
The Lie symmetry is in this way conserved in the isolated systems but
the composed ones must be considered as a whole: sub-systems
cannot be secluded and this basis could be an algebraic way to introduce the
interaction.

Let us call $\{I_3, I_\pm\}$  the generators inside the algebra
$U_q(su(2))$ that close the Lie commutation relations
(\ref{commL}).
Of course $\{I_m\}$ do not close a Lie-Hopf algebra also if
the all Lie algebraic relations
(\ref{commL}, \ref{reprL}, \ref{CasimirL}) remain valid.
Indeed these relations are 
consistent with the primitive co-product (\ref{DL}) but
also with $\infty$-many others. In our case the reference 
structure is $U_q(su(2))$ and the coalgebra imposes

\begin{equation}\label{DI}
\begin{array}{l}
\Delta(I_3) = 1 \otimes I_3 + I_3 \otimes 1 + o(3)\\[0.3cm]
\Delta(I_\pm) =  (1 \otimes I_\pm + I_\pm \otimes 1) + \,z\, 
(I_3 \otimes I_\pm - I_\pm \otimes I_3) + o(3) \;.
\end{array}\end{equation}
Formulas (\ref{mapKtoJ}) are rewritten for $z'=0$ as

\[
I_3 = J_3 \qquad
I_{\pm} = 
 \sqrt{[1]_q} \;
\sqrt{\frac{ (j+1/2)^2 - (J_3\mp 1/2)^2}{[j+1/2]_q^{\,2} 
- [J_3\mp 1/2]_q^{\,2}}} \;
J_{\pm}
\]
and perturbative formulas are
\[
I_\pm = J_\pm - 
\frac{z^2}{3} {\cal S}(J_3^2 J_\pm)- 
\frac{z^2}{6} {\cal S}(J_\pm^2 J_\mp) + o(4) \;, 
\]
\[
J_\pm = I_\pm + 
\frac{z^2}{3} {\cal S}(I_3^2 I_\pm)+
\frac{z^2}{6} {\cal S}(I_\pm^2 I_\mp) + o(4) \,,
\]
\[
\Delta(I_\pm) = \Delta(J_\pm)
 - 
\frac{z^2}{3} \Delta\bigl({\cal S}(J_3^2 J_\pm)\bigr)- 
\frac{z^2}{6} \Delta\bigl({\cal S}(J_\pm^2 J_\mp)\bigr)   
 + o(4) \,,
\]
\noindent such that

\begin{equation*}
\begin{array}{l}
\Delta(I_\pm) = q^{I_3} \otimes I_\pm + 
I_\pm \otimes q^{-I_3}  \;-  \\[0.3cm]
(z^2/3) \Bigl[\Delta_{(0)}({\cal S}(I_3^2 I_\pm))- 
1 \otimes {\cal S}(I_3^2 I_\pm)- {\cal S}(I_3^2 I_\pm)\otimes 1 
\Bigr] \;- \\[0.3cm]
(z^2/6) \Bigl[\Delta_{(0)}({\cal S}(I_\pm^2 I_\mp))- 
1 \otimes {\cal S}(I_\pm^2 I_\mp)- {\cal S}(I_\pm^2 I_\mp)\otimes 1 
\Bigr]
 + o(4)\; ,
\end{array}\end{equation*}
in agreement with eqs.(\ref{DI}):
in spite of the Lie-like commutation relations of the $\{I_m\}$, the
Hopf algebra is always $U_q(su(2))$ and not $U(su(2))$.

{\section{Crystal Basis $su_{\infty}(2)$}

The crystal basis of $U_q(su(2))$ can be considered as the conjugate of 
the Lie one: it is the quantum basis for $z'=\infty$ in the 
Riemann sphere. Of course, 
because of the singularities of polynomials and
exponentials at the point $z'=\infty$,
we have to renormalize the rising and lowering operators

\[
C_3^{\;q'} := K_3 \qquad\quad C_\pm^{\;q'} :=  
\frac{\sqrt{[1]_{q'}}}{[j+1/2]_{q'}} \;K_\pm \;,
\]
such that
\[
C_3^{\;q'} \,|j m\rangle = |j m\rangle\, m \qquad\;
C_\pm^{\, q'} \;|j m\rangle\, =\; |j ~m\pm 1\rangle \;\;\sqrt{1 - 
\frac{[m\pm1/2]_{q'}^{\;2}}{[j+1/2]_{q'}^{\;2}}}\, .
\]
The crystal basis is then defined as:
\[
C_3:=C_3^{\;q'} \quad\qquad C_\pm := \lim_{z' \to \infty} C_\pm^{\, q'}\,. 
\]
The representations of $\{C_m\}$ are thus 
\[
C_3|j,m\rangle = |j, m\rangle ~m \qquad\;\;
C_\pm |j,m\rangle = |j, m\pm 1\rangle ~~(1-\delta_{j,\pm m})
\]
and the mapping on $\{J_m\}$ is
\be\label{mapCtoJ}
C_\pm = \frac{\sqrt{[1]_q}}
{\sqrt{[j+1/2]_q^{\;2} - 
[J_3\mp 1/2]_q^{\;2}}}  \;J_\pm
\ee
where a consistent definition has been adopted on the highest weight states.

The co-product is now calculated as

\[
\Delta(C_\pm) =
\frac{\sqrt{[1]_q}}
{\sqrt{\Delta([j+1/2]_q^{\;2}) - 
\Delta([J_3\mp 1/2]_q^{\;2})}} \;\, \Delta(J_\pm)
\]
where, as before, $\Delta$ in the {\it rhs} is determined from 
eqs. (\ref{coprodJ}) and $\Delta(C_\pm)$ is calculated in terms of
$\{\Delta(J_m)\}$.
Inverting eq. (\ref{mapCtoJ}) we have
\[
J_\pm = \sqrt{[j+1/2]_q^2 - 
[C_3\mp 1/2]_q^2} \;\; C_\pm \,/\, \sqrt{[1]_q}\;.
\]
Because the Casimir can be written as
\[
j = \sum_{k=0}^\infty (C_-)^k C_3 (1- C_- C_+) (C_+)^k \qquad \in U_q(su(2))\;,
\]
\noindent $\Delta(C_\pm)$ can be written in function of $\{C_m \otimes 1\}$
and $\{1 \otimes C_m\}$\,.




\end{document}